\documentclass[12pt,a4paper]{article}
\usepackage{mathrsfs}
\usepackage{epsfig, graphicx}
\usepackage{latexsym,amsfonts,amsbsy,amssymb}
\usepackage{amsmath,amsthm}
\usepackage{hyperref}
\usepackage{color}
\makeatletter 

\@addtoreset{equation}{section}
\makeatother
\allowdisplaybreaks 

\newtheorem{lemma}{Lemma}[section]

\newtheorem{algorithm}{Algorithm}[section]

\begin{document}
\title{Convergence Analysis of Shift-Inverse Method with Richardson Iteration For Eigenvalue Problem\thanks{This work is
supported in part by National Natural Science Foundations
of China (NSFC 91330202, 11001259, 11371026, 11201501, 11031006, 2011CB309703, 11171251) and
the National Center for Mathematics and Interdisciplinary Science, CAS.}}
\author{
Yunhui He\thanks{Department of Mathematics and Statistics,
Memorial University of Newfoundland, St. John's, NL A1C 5S7, Canada ({\tt  yunhui.he@mun.ca})}
\ \ \ and\ \
Hehu Xie\thanks{LSEC, NCMIS, Academy of Mathematics and Systems Science,
Chinese Academy of Sciences, Beijing 100190, China
({\tt hhxie@lsec.cc.ac.cn})}
}
\date{}
\maketitle
\begin{abstract}
In this paper, we consider the shift-inverse method with Richardson iteration step for the eigenvalue problems. It will
be shown that the convergence speed depends heavily on the eigenvalue gap between the desired eigenvalue and undesired ones.

\end{abstract}
\section{Introduction}
This paper is to discuss the convergence behavior of the Richardson iteration step in the shift-inverse (or Rayleigh quotient) method
for the eigenvalue problem. It is well known that the shift-inverse method has a superlinear convergence order if we have
good enough initial guesses (cf. \cite{ChenHeLiXie,Notay,Saad}). But in the shift-inverse method, we need to solve the almost singular
linear system which is not so easy work. Thus it is necessary to consider the effectiveness of the iteration steps for the
almost singular systems.
The first aim of this paper is to discuss the Richardson iteration step applied to the linear equations
deduced from the shift-inverse method for the eigenvalue problem. The common idea is that it should be difficult
to solve the linear equations since they are almost singular when the shift close to a singular.
We will show here that the common iteration method can also work very well for these linear equations
since the final aim is to find the eigenvectors rather than solving the deduced linear equations.
The most important property is that the iteration step can reduce the relative scale of the undesired eigenvectors
according to the desired eigenvectors. After normalization, the errors in the undersired eigenvectors
will become smaller and smaller.

The rest of this paper is organized as follows. In the next section, we introduce the shift-inverse method and
the Richardson iteration for the shift-inverse equation will be analyzed in Section 3.
Section 4 is devoted to a numerical example and some concluding remarks will be given in the final section.
\section{Shift-inverse method}
For simplicity, we solve the following standard eigenvalue problem:
Find $(\lambda,x)\in \mathbb{R}\times \mathbb{R}^n$ such that
\begin{eqnarray}\label{Eigenvalue_Problem_Standard}
Ax&=&\lambda x.
\end{eqnarray}
The shift-inverse method for the eigenvalue problem can be descried as follows:

\begin{algorithm}\label{Algorithm_Shift_Inverse}
Assume we have the eigenpair approximations
$(\lambda_i^{(k)},x_i^{(k)})$  $(i=1,\cdots,\ell)$.
\begin{enumerate}
\item For $i=1,\cdots,\ell$, solve the following linear equation
\begin{eqnarray}\label{Linear_Equation_Shift_Inverse}
(A-\tau_kI)\tilde{x}_i^{(k)}&=&x_i^{(k)}.
\end{eqnarray}
\item
Build the following small scale eigenvalue problem
\begin{eqnarray}
(X^{(k)})^{T}AX^{(k)}z_i&=&\lambda_i^{(k+1)}(X^{(k)})^{T}X^{(k)}z_i,\ \ \ \
i=1,\cdots,\ell,
\end{eqnarray}
where $X^{(k)}=[\tilde{x}_1^{(k)},\cdots,\tilde x_{\ell}^{(k)}]$. Solve this eigenvalue problem
to obtain the new eigenpair approximations $(\lambda_i^{(k+1)},x_i^{(k+1)})=(\lambda_i^{(k+1)},X^{(k)}z_i)$
 $(i=1,\cdots,\ell)$.
\end{enumerate}
\end{algorithm}
It is well known that the shift-inverse iteration is a basic numerical method for the
eigenvalue problem.
\begin{lemma}\label{Convergence_Lemma}(\cite[Theorem 1 and Remark 1]{ChenHeLiXie})
The shift-inverse method defined by Algorithm \ref{Algorithm_Shift_Inverse} converges
linearly. Especially, if we we choose $\tau_k=\lambda_1^{(k)}$, the shift-inverse
method has the cubic convergence
\begin{eqnarray}\label{Convergence_Rate}
\|x-x^{(k+1)}\|_A &\leq & C\|x-x^{(k)}\|_A^3.
\end{eqnarray}
\end{lemma}

\section{Richardson iteration for the shift-inverse equation}
In this section, we give a detailed analysis for the solution of the linear system (\ref{Linear_Equation_Shift_Inverse})
with the basic Richardson iteration.

Assume that the eigenvalues of $A$ satisfy
\begin{eqnarray*}
\lambda_{\min}(A)=\lambda_{1}\leq \lambda_{2}\leq\cdots \leq \lambda_{\ell}< \lambda_{\ell+1}\leq \cdots\leq\lambda_{n}
\end{eqnarray*}
and corresponding eigenvectors
\begin{eqnarray*}
\phi_{1},\ \phi_{2},\ \cdots,\ \phi_{n}
\end{eqnarray*}
with $\phi_i^T\phi_j=\delta_{ij}$ and $\delta_{ij}$ denotes the Kronecker function.
It means that the vector system $\{\phi_1, \cdots, \phi_n\}$ is an orthonormal basis for the
linear space $\mathcal{R}^n$.

Let the vector $x_i^{(k)}$ has the following expansion form
\begin{eqnarray}\label{Decomposition_X_k}
x_i^{(k)}&=&\sum_{j=1}^n\alpha_{ij}^{(k)}\phi_j,\ \ \ \ i=1,\cdots,\ell,
\end{eqnarray}
and $\alpha_{ij}^{(k)}$ ($j=1,\cdots,n$) satisfy the normalization condition
\begin{eqnarray}\label{Normalization_Condition}
\sum_{j=1}^n(\alpha_{ij}^{(k)})^2&=&1,\ \ \ i=1,\cdots,\ell.
\end{eqnarray}

Different from the iteration method for the linear equation, the aim of the iteration step
for the linear system (\ref{Linear_Equation_Shift_Inverse}) is to reduce the terms
$\alpha_{ij}^{(k)}$ with $j=\ell+1,\cdots,n$ corresponding to the undesired eigenvalues
$\lambda_{\ell+1},...,\lambda_n$. The idea here is to reduce the relative scales of all the terms
$\alpha_{ij}^{(k)}$ with $j=\ell+1,\cdots,n$ according to the desired
terms $\alpha_{ij}^{(k)}$ with $j=1,\cdots,\ell$.

For simplicity, we consider the following Richardson iteration scheme:
\begin{eqnarray}\label{Richardson_Shift}
 x_i^{(k+1)}&=&x_i^{(k)}-\theta\big((A-\tau_k I)x_i^{(k)}-x_i^{(k)}\big)\nonumber\\
 &=&\Big(\big(1+\theta(1+\tau_k)\big)I-\theta A\Big)x_i^{(k)},
\end{eqnarray}
where $\theta\in (0,1)$.
The iteration matrix is
\begin{eqnarray}
G_k&=&\big(1+\theta(1+\tau_k)\big)I-\theta A.
\end{eqnarray}
The iteration scheme (\ref{Richardson_Shift}) has the following form
\begin{eqnarray}
x_i^{(k+1)}&=&G_kx_i^{(k)}=\sum_{j=1}^n\alpha_{ij}^{(k)}G_k\phi_j\nonumber\\
&=&\sum_{j=1}^n\big(1+\theta(1+\tau_k)-\theta\lambda_j\big)\alpha_{ij}^{(k)}\phi_j.
\end{eqnarray}
It means that the following relations hold for $i=1,\cdots,\ell$
\begin{eqnarray}\label{Relation_k_k+1}
\alpha_{ij}^{(k+1)}=\big(1+\theta(1+\tau_k)-\theta\lambda_j\big)\alpha_{ij}^{(k)},
\ \ \ \ \ \  j=1,\cdots,n.
\end{eqnarray}
Based on the relation (\ref{Relation_k_k+1}), it is reasonable
that the convergence rate can be estimated as follows
\begin{eqnarray}
{\rm Rate}_{k+1}:=\frac{\max_{\lambda_j\geq\lambda_{\ell+1}}\big|1+\theta(1+\tau_k)-\theta\lambda_j\big|}
{\min_{\lambda_j\leq \lambda_{\ell}}\big|1+\theta(1+\tau_k)-\theta\lambda_j\big|}
\end{eqnarray}
In order to minimize the value ${\rm Rate}_{k+1}$, we chose $\tau_k$ such that $|1+\theta(1+\tau_k)-\theta\lambda_{\ell+1}|=|1+\theta(1+\tau_k)-\theta\lambda_n|$.
It means that we have
\begin{eqnarray}
\tau_k=\frac{\lambda_{\ell+1}+\lambda_n}{2}-\frac{1}{\theta}-1
\end{eqnarray}
and
\begin{eqnarray}\label{Convergence_Rate_Result}
{\rm Rate}_{k+1} = \frac{\lambda_n-\lambda_{\ell+1}}{\lambda_n+\lambda_{\ell+1}-2\lambda_{\ell}}
=1-2\frac{\lambda_{\ell+1}-\lambda_{\ell}}{\lambda_n+\lambda_{\ell+1}-2\lambda_{\ell}}.
\end{eqnarray}
The convergence rate (\ref{Convergence_Rate_Result})
shows that the convergence rate depends strongly on the eigenvalue gap $\lambda_{\ell+1}-\lambda_{\ell}$
according to the value $\lambda_n+\lambda_{\ell+1}-2\lambda_{\ell}$.

\section{Numerical results}
In this section, we present a simple example to illustrate the relation between the convergence speed
and the eigenvalue gap which has been given in (\ref{Convergence_Rate_Result}). Here we solve the eigenvalue
problem: Find $(\lambda,x)\in \mathbb R\times \mathbb R^4$ such that
\begin{eqnarray*}
Ax=\lambda x,
\end{eqnarray*}
where $A={\rm diag}(1,2,2.01,4)$.

\begin{figure}[ht]
\centering
\includegraphics[width=7cm,height=7cm]{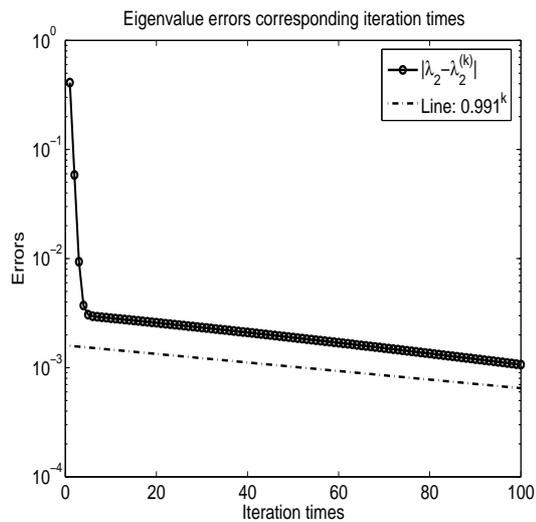}
\caption{The eigenvalue errors of the shift-inverse method with Richardson iteration step}\label{Error_Figure}
\end{figure}

We choose $\theta=0.5$, $\tau = \frac{2.01+4}{2}-\frac{1}{0.5}-1=0.005$ and solve the eigenvalue problem to
obtain the first two eigenpair approximations according to the eigenvalues $\lambda_1=1$ and $\lambda_2=2$.
We choose two random vectors as the initial eigenvectors to do the shift-inverse method with the Richardson
iteration step described in Algorithm \ref{Algorithm_Shift_Inverse}. The error of the second eigenvalue approximations
are presented in Figure \ref{Error_Figure} which shows the slow convergence speed since the gap $\lambda_3-\lambda_2$
is very small.

\section{Concluding remarks}
In this paper, we discuss the convergence behavior of the shift-inverse method with Richardson
iteration step for the eigenvalue problem. It is shown that the eigenvalue gap
decide the convergence speed.


\bibliographystyle{amsplain}

\end{document}